\documentclass[12pt]{article}

\usepackage{amsmath, epsfig, cite}
\usepackage{amssymb}
\usepackage{amsfonts}
\usepackage{latexsym}
\usepackage{amsthm}

\makeatletter
\renewcommand{\@seccntformat}[1]{{\csname the#1\endcsname}.\hspace{.5em}}
\makeatother

\newtheorem{thm}{Theorem}[section]

\newtheorem{cor}[thm]{Corollary}

\newtheorem{lem}[thm]{Lemma}

\newcommand{\pf}{\noindent{\it Proof.} }

\renewcommand{\qed}{\hfill$\Box$\medskip}

\numberwithin{equation}{section}

\begin{document}

\renewcommand{\thefootnote}{*}

\begin{center}
{\Large\bf A $q$-Clausen-Orr type formula and its applications}
\end{center}

\vskip 2mm \centerline{Hong-Fang Guo$^{1}$, Victor J. W. Guo$^2$\footnote{Corresponding author.}  and Jiang Zeng$^{3}$}
\begin{center}
{\footnotesize $^1$Department of Mathematics, East China Normal University, Shanghai
200241,\\
People's Republic of China\\[5pt]

$^2$School of Mathematical Sciences, Huaiyin Normal University, Huai'an, Jiangsu 223300,\\
 People's Republic of China\\[5pt]

$^3$Univ Lyon, Universit\'e Claude Bernard Lyon 1,  CNRS UMR 5208,   Institut Camille
Jordan,  43, blvd. du 11 novembre 1918,
F-69622 Villeurbanne Cedex, France\\[5pt]
{\tt ghfkeji@126.com,     jwguo@hytc.edu.cn,    zeng@math.univ-lyon1.fr}}
\end{center}

\vskip 0.7cm \noindent{\small{\bf Abstract.}}
We show  that certain terminating
$_{6}\phi_5$ series can be factorized into a product of two $_{3}\phi_{2}$ series. As applications we  prove
a summation formula for a product of two $q$-Delannoy numbers along with
 some congruences for sums involving $q$-Delannoy numbers.
 This confirms three recent conjectures of the second author.


\section{Introduction}
Clausen's formula
\begin{align*}
&\left\{_{2}F_{1}\left[\begin{array}{c}
a, b\\
a+b+1/2
\end{array};z\right]
\right\}^2
={}_3F_{2}\left[\begin{array}{c}
2a, 2b, a+b\\
a+b+1/2, 2a+2b
\end{array};z
\right]
\end{align*}
plays a central role in Ramanujan's derivation for various series for $1/\pi$.
See \cite{AVZ11,CTYZ11} for some recent developments of this formula.
More general formulas connecting products of
two hypergeometric series as a single series were obtained by Orr in 1899 (see \cite[p.~75]{Slater}).
%
This paper was motivated by
 a recent paper of  the second author~\cite{Guo}, where he  proved some
 congruences of sums involving even powers of Delannoy numbers
and raised some problems of finding the $q$-analogues.  Recall that
the Delannoy numbers count lattice paths from ($0,0$) to ($n,m$) consisting of horizontal (1,0), vertical (0,1), and diagonal (1,1) steps,  and have the following explicit formulas in terms of binomial coefficients:
\begin{align}
D(m,n):=\sum_{k=0}^{n}{n\choose  k}{n+m-k\choose n}=\sum_{k=0}^{n}{n\choose  k}{m\choose k}2^k. \label{delannoynumber}
\end{align}
The reader is referred to Dziemia\'nczuk~\cite{Dzie} and the references therein  for
how to generalize Delannoy numbers via counting weighted lattice paths.
Recall that the {\it basic hypergeometric series $_{r}\phi_s$}  is defined as
$$
_{r}\phi_{s}\left[\begin{array}{c}
a_1,a_2,\ldots,a_{r}\\
b_1,\ldots,b_{s}
\end{array};q, z
\right]
=\sum_{k=0}^{\infty}\frac{(a_1;q)_k(a_2;q)_k\cdots (a_r;q)_k}
{(q;q)_k (b_1;q)_k \cdots (b_s;q)_k} \left((-1)^k q^{k\choose 2}\right)^{1+s-r} z^k,
$$
where  $(a;q)_n=(1-a)(1-aq)\cdots(1-aq^{n-1})$  for $n=1,2,\ldots,$  and $(a;q)_0=1$.
Aiming to answer the $q$-problems in \cite{Guo}, we  are led to prove the following $q$-Clausen-Orr type formula.

\begin{thm}\label{thm:1.4}
Let $n$ be a non-negative integer. Then
\begin{align}
_{3}\phi_{2}\left[\begin{array}{c}
q^{-n}, a, x\\
c,\, 0
\end{array};q,q
\right]
{}_{3}\phi_{2}\left[\begin{array}{c}
q^{-n}, a, c/x\\
c,\, 0
\end{array};q,q
\right]
=a^n {}_{6}\phi_{5}\left[\begin{array}{c}
q^{-n}, cq^{n}, a, c/a, x, c/x\\
c,\sqrt{c}, -\sqrt{c}, \sqrt{cq}, -\sqrt{cq}
\end{array};q,q
\right].  \label{eq:q-new-Clausen}
\end{align}
\end{thm}

It is interesting to compare \eqref{eq:q-new-Clausen} with
Jackson's $q$-analogue of Clausen's formula \cite{Jackson40,Jackson41}:
\begin{align*}
&_{2}\phi_{1}\left[\begin{array}{c}
a, b\\
abq^{\frac{1}{2}}
\end{array};q,z
\right]
{}_{2}\phi_{1}\left[\begin{array}{c}
a, b\\
abq^{\frac{1}{2}}
\end{array};q,zq^{\frac{1}{2}}
\right]
={}_{4}\phi_{3}\left[\begin{array}{c}
a, b, a^{\frac{1}{2}}b^{\frac{1}{2}}, -a^{\frac{1}{2}}b^{\frac{1}{2}}\\
ab, a^{\frac{1}{2}}b^{\frac{1}{2}}q^{\frac{1}{4}}, -a^{\frac{1}{2}}b^{\frac{1}{2}}q^{\frac{1}{4}}
\end{array};q^{\frac{1}{2}},z
\right],\ |z|<1,
\end{align*}
and the following  $q$-analogue of Clausen's formula \cite[Appendix (III.22)]{GR}:
\begin{align}
\left\{{}_{4}\phi_{3}\left[\begin{array}{c}
a, b, abz, ab/z\\
abq^{\frac{1}{2}}, -abq^{\frac{1}{2}}, -ab
\end{array};q,q
\right]\right\}^2
={}_{5}\phi_{4}\left[\begin{array}{c}
a^2, b^2, ab, abz, ab/z\\
abq^{\frac{1}{2}}, -abq^{\frac{1}{2}}, -ab, a^2b^2
\end{array};q,q
\right], \label{eq:classic-qClausen}
\end{align}
where both  series  are supposed to be terminated.
Indeed, letting $c=x^2$, the identity \eqref{eq:q-new-Clausen} reduces to
the following formula, which seems to be new.
\begin{cor} Let $n$ be a non-negative integer. Then
\begin{align}
\left\{_{3}\phi_{2}\left[\begin{array}{c}
q^{-n}, a, x\\
x^2,\, 0
\end{array};q,q
\right]\right\}^2
=a^n {}_{5}\phi_{4}\left[\begin{array}{c}
q^{-n}, x^2q^{n}, a, x^2/a, x\\
x^2,-x, xq^{\frac{1}{2}}, -xq^{\frac{1}{2}}
\end{array};q,q
\right].  \label{eq:final-square}
\end{align}
\end{cor}
In particular,  the right-hand side of \eqref{eq:final-square} is non-negative for real $a$, $x$, and $q$.
Furthermore,
if $n$ is even, then by \eqref{eq:classic-qClausen}, the right-hand side of \eqref{eq:final-square} may be written as
\begin{align*}
a^n\left\{{}_{4}\phi_{3}\left[\begin{array}{c}
q^{-\frac{n}{2}}, xq^{\frac{n}{2}}, a, x^2/a\\
xq^{\frac{1}{2}}, -xq^{\frac{1}{2}}, -x
\end{array};q,q
\right]\right\}^2.
\end{align*}
Writing $n=2m$ in \eqref{eq:final-square}
and taking the square root  we obtain
an identity between two  polynomials in $a$ of degree $2m$, where the sign is determined by
comparing the coefficients of $a^{2m}$. We record  the resulting  formula  as the second  corollary.
\begin{cor} Let $m$ be a non-negative integer. Then
\begin{align*}
_{3}\phi_{2}\left[\begin{array}{c}
q^{-2m}, a, x\\
x^2,\, 0
\end{array};q,q
\right]
=a^m
{}_{4}\phi_{3}\left[\begin{array}{c}
q^{-m}, xq^{m}, a, x^2/a\\
xq^{\frac{1}{2}}, -xq^{\frac{1}{2}}, -x
\end{array};q,q
\right].
\end{align*}
\end{cor}
For some other $q$-Clausen type formulas,  the reader is referred to
Gasper and Rahman\cite[Exercise 8.17]{GR} and  Schlosser~\cite{Schlosser}.
On the other hand, in their study of some $q$-supercongruences for certain truncated basic hypergeometric series related to \cite{SunZH2,Hamme},
Guo and Zeng~\cite{GZ2}
stumbled  on the following $q$-Clausen-Orr  type formula:
\begin{align}
&\hskip -3mm \left(\sum_{k=s}^{n} \frac{(q^{-2n};q^2)_k (x;q)_k q^k }{(q;q)_{k-s} (q;q)_{k+s}}\right)
  \left(\sum_{k=s}^{n} \frac{(q^{-2n};q^2)_k (q/x;q)_k q^k }{(q;q)_{k-s} (q;q)_{k+s}}\right) \nonumber \\
&=\frac{(-1)^n(q^2;q^2)_n^2 q^{-n^2}}{(q^2;q^2)_{n-s}(q^2;q^2)_{n+s}}
\sum_{k=s}^{n} \frac{(-1)^k(q^2;q^2)_{n+k} (x;q)_k (q/x;q)_k  q^{k^2-2nk}}{(q^2;q^2)_{n-k} (q;q)_{k-s}(q;q)_{k+s} (q;q)_{2k}}.
\label{eq:important-thm-0}
\end{align}
Noticing that $(c;q)_{2k}=(\sqrt{c};q)_k(-\sqrt{c};q)_k(\sqrt{cq};q)_k(-\sqrt{cq};q)_k$,
we can rewrite \eqref{eq:q-new-Clausen} as
\begin{align}
&\hskip -3mm \left(\sum_{k=0}^{n} \frac{(q^{-n};q)_k (a;q)_k(x;q)_k q^k }{(q;q)_{k} (c;q)_{k}}\right)
  \left(\sum_{k=0}^{n} \frac{(q^{-n};q)_k (a;q)_k(c/x;q)_k q^k }{(q;q)_{k} (c;q)_{k}}\right) \nonumber \\
&=a^n
\sum_{k=0}^{n} \frac{(q^{-n};q)_k (cq^{n};q)_k (a;q)_k (c/a;q)_k (x;q)_k (c/x;q)_k  q^{k}}{(q;q)_{k}(c;q)_{k} (c;q)_{2k}}. \label{eq:thm1.4}
\end{align}
Clearly \eqref{eq:thm1.4} is an extension of \eqref{eq:important-thm-0}.
In the same vein we shall establish the following result.
\begin{thm}\label{thm:last}Let $n$ be a non-negative integer. Then
\begin{align}
&\hskip -3mm \left(\sum_{k=0}^{n} \frac{(q^{-n};q)_k (x;q^2)_k q^k }{(q;q)_{k} (c;q)_{k}}\right)
 \left(\sum_{k=0}^{n} \frac{(q^{-n};q)_k (x;q^2)_k c^k q^{nk-{k\choose 2}} }{(q;q)_{k} (c;q)_{k} x^k }\right) \nonumber \\
&=\sum_{k=0}^{n} \frac{(q^{-n};q)_k (cq^{n};q)_k (x;q^2)_k (c^2/x;q^2)_k  q^{k}}{(q;q)_{k}(c;q)_{k} (c;q)_{2k}}. \label{eq:special-3}
\end{align}
\end{thm}

In this paper we shall consider two $q$-analogues of $D(m,n)$.  We first recall some standard $q$-notation (see \cite{GR}).
The {\it $q$-binomial coefficients} are given by
\begin{align*}
{n\brack k}:=
\begin{cases}
\displaystyle\frac{(q;q)_n}{(q;q)_k(q;q)_{n-k}}
&\text{if $n\geqslant k\geqslant 0$,} \\[10pt]
0 &\text{otherwise.}
\end{cases}
\end{align*}
The following two  natural $q$-analogues of Delannoy numbers were introduced in
\cite[p.~30]{Dzie} and \cite{Pan15}:
\begin{align}
D_q(m,n)&:=\sum_{k=0}^{n}q^{k\choose 2}{n\brack k}{n+m-k\brack n},\\
D^*_q(m,n)&:=\sum_{k=0}^{n}q^{k+1\choose 2}{n\brack k}{n+m-k\brack n}.
\end{align}
Note that $D^*_{q}(m,n):=q^{mn}D_{q^{-1}}(m,n)$.
We first show that  both $D_{q}(m,n)$ and $D^*_{q}(m,n)$ have a $q$-analogue of
the second expression in \eqref{delannoynumber} and provide a $q$-analogue of
 \cite[(3.1)]{Guo}, which was asked in  \cite[Problem 5.2]{Guo}.
\begin{thm}\label{thm:1}
Let $m$ and $n$ be non-negative integers. Then
\begin{align}
D_q(m,n)&=\sum_{k=0}^{m}q^{(m-k)(n-k)}{m\brack k}{n\brack k}(-1;q)_k,\label{D1}\\
D^*_q(n,m)&=\sum_{k=0}^{m}q^{(m-k)(n-k)}{m\brack k}{n\brack k}(-q;q)_k.\label{D2}
\end{align}
Moreover,
\begin{align}
D_q(m,n)D^*_{q}(m,n)=\sum_{k=0}^{n}q^{(m-k)(n-k)}{n+k\brack 2k}{m\brack k}{m+k\brack k}(-1;q)_k(-q;q)_k. \label{eq:thm:1}
\end{align}
\end{thm}

Applying the formula \eqref{eq:thm:1}, we shall prove the following two results originally conjectured by the second author \cite[Conjectures 5.3 and 5.4]{Guo}.
\begin{thm}\label{thm:2}
Let $p$ be an odd prime and $m$ a positive integer. Then
\begin{align}
&\hskip -3mm \sum_{k=0}^{p-1}\frac{1-q^{2k+1}}{1-q}D_q(m,k) D_{q^{-1}}(m,k) q^{-k} \nonumber\\
&\equiv\begin{cases}
\displaystyle\frac{1-q^{-2m}}{1-q^2}q\pmod{[p]^2}&\text{if $m\equiv 0\pmod{p}$}, \\[10pt]
\displaystyle\frac{1-q^{2m+2}}{1-q^2}q\pmod{[p]^2}&\text{if $m\equiv -1\pmod{p}$}, \\[10pt]
0\pmod{[p]^2}&\text{otherwise},
\end{cases}\label{eq:cases}
\end{align}
where $[p]=1+q+\cdots+q^{p-1}$, and the congruences are understood in the polynomial ring $\mathbb{Z}[q]$.
\end{thm}

\begin{thm}\label{thm:3}
Let $m$, $n$, and $r$ be positive integers. Then all of
\begin{align}
&\sum_{k=0}^{n-1}\frac{(1-q^{m})(1-q^{m+1})(1-q^{2k+1})}{(1-q^2)(1-q^n)^2}D_q(m,k)D_{q^{-1}}(m,k)q^{-k}, \label{eq:thm3-1}\\[5pt]
&\sum_{k=0}^{n-1}\frac{1-q^{2k+1}}{1-q^n}D_q(m,k)^r D_{q^{-1}}(m,k)^r q^{-k},  \label{eq:thm3-2}\\[5pt]
&\sum_{k=0}^{n-1}(-1)^{n-k-1}\frac{1-q^{2k+1}}{1-q^n}D_q(m,k)^r D_{q^{-1}}(m,k)^r q^{k\choose 2}  \label{eq:thm3-3}
&\end{align}
are Laurent polynomials in $q$ with non-negative integer coefficients.
\end{thm}
Note  that Theorem \ref{thm:3} is a $q$-analogue of \cite[Theorem 1.1]{Guo} for the first three polynomials.
The rest of the paper is organized as follows.
We shall give three lemmas in Section~\ref{sec:lemmas} and prove Theorem~\ref{thm:1.4} in
Section~\ref{sec:thm:1.4}. In Sections~\ref{Pr:thm:last} and \ref{sec:proofoftheorem1} we prove
Theorems~\ref{thm:last} and \ref{thm:1}.
In Sections~\ref{sec:pf2} and \ref{sec:pf3}, by using Theorem~\ref{thm:1} we give proofs of
Theorems~\ref{thm:2} and \ref{thm:3} respectively.

\section{Three lemmas\label{sec:lemmas}}
The following three lemmas are crucial ingredients of our proof of Theorem {\rm\ref{thm:1.4}.
\begin{lem}Let $n$ and $h$ be positive integers and let $m$ be a non-negative integer with $h\leqslant n-m$. Then
\begin{align}
&\sum_{j=0}^{m}\sum_{k=0}^n \frac{(q^{-n};q)_{j}  (q^{-n};q)_{k}(x;q)_{j} (x;q)_{k} (q^{j-m-h+1};q)_{h-1}(q^{k-m-h+1};q)_{h-1} (1-q^{k-j}) q^{2j+k}  }
{(q;q)_{j} (q;q)_{k} (c;q)_{j} (c;q)_{k}} \notag\\
&=\frac{(q;q)_n (q;q)_{h-1} (x;q)_{m+h}(c/x;q)_{n-h}x^{n-h}q^{\frac{m^2+3m}{2}-mn-mh-h^2+h}}
{(-1)^{m-1}(q;q)_{m} (c;q)_{m} (c;q)_{n}(q;q)_{n-m-h}}. \label{eq:last-lem}
\end{align}
\end{lem}

\pf  Note that both sides of \eqref{eq:last-lem} are polynomials in $x$ of degree $m+n$ with
the same leading coefficient. Therefore, to prove \eqref{eq:last-lem}, it
suffices to prove that both sides have the same roots as polynomials in $x$.
Denote the left-hand side of \eqref{eq:last-lem} by $L_{m,n}(x)$. We first assert that
\begin{align}
L_{m,n}(x)&=\sum_{j=0}^{m}\sum_{k=m+1}^n \frac{(q^{-n};q)_{j}  (q^{-n};q)_{k}(x;q)_{j} (x;q)_{k}  }
{(q;q)_{j} (q;q)_{k}(c;q)_{j} (c;q)_{k}} \notag\\
&\quad\times (q^{j-m-h+1};q)_{h-1}(q^{k-m-h+1};q)_{h-1} (1-q^{k-j}) q^{2j+k}. \label{eq:last-lem-2}
\end{align}
In fact, since $(1-q^{j-k})q^{2k+j}=-(1-q^{k-j})q^{k+2j}$, we have $\sum_{j=0}^{m}\sum_{k=0}^{m}=0$
for the summands in $L_{m,n}(x)$.
We now consider the following two cases.
\begin{itemize}

\item For $x=q^{-r}$ with $0\leqslant r\leqslant m+h-1$, we have
\begin{align*}
L_{m,n}(q^{-r})&=\sum_{j=0}^{m}\sum_{k=0}^n \frac{(q^{-n};q)_{j}  (q^{-n};q)_{k}(q^{-r};q)_{j} (q^{-r};q)_{k}   }
{(q;q)_{j} (q;q)_{k}(c;q)_{j} (c;q)_{k}} \notag\\
&\quad{}\times (q^{j-m-h+1};q)_{h-1}(q^{k-m-h+1};q)_{h-1} (1-q^{k-j}) q^{2j+k}.
\end{align*}
If $r\leqslant m$, then $L_{m,n}(q^{-r})=0$
by the antisymmetry of $j$ and $k$ in $L_{m,n}(q^{-r})$. If $r\geqslant m+1$, then $h\geqslant r-m+1$,
i.e., $r-m-h+1\leqslant 0$, and so $(q^{k-m-h+1};q)_{h-1}=0$ for $m+1\leqslant k\leqslant r$. Hence, by \eqref{eq:last-lem-2},
we again get $L_{m,n}(q^{-r})=0$.

\item For $x=cq^{r}$ with $0\leqslant r\leqslant n-h-1$, we shall prove that
\begin{align}
\sum_{k=0}^n \frac{(q^{-n};q)_{k}(cq^r;q)_{k} (q^{k-m-h+1};q)_{h-1} (1-q^{k-j}) q^{2j+k}  }
{(q;q)_{k} (c;q)_{k}} =0.  \label{eq:last-lem-3}
\end{align}
In fact, we can rewrite the left-hand side of \eqref{eq:last-lem-3} as
\begin{align}
\sum_{k=0}^{n} (-1)^{k}{n\brack k}q^{-nk+{k\choose 2}} R_k,\label{eq:last-lem-4}
\end{align}
where
$$
R_k=
\frac{(cq^{r};q)_{k} (q^{k-m-h+1};q)_{h-1} (1-q^{k-j}) q^{2j+k}} {(c;q)_{k}}.
$$
Since
\begin{align*}
\frac{(cq^{r};q)_{k}}{(c;q)_{k}}=
\frac{(cq^{k};q)_{r}} {(c;q)_{r}},
\end{align*}
we see that $R_k$ is a polynomial in $q^k$ of degree $r+h-1+2\leqslant n$ with no constant term.
By the $q$-binomial theorem (see, for example, \cite[Theorem 3.3]{Andrews})
\begin{align}
\sum_{k=0}^{n}(-1)^{k}{n\brack k}q^{{k\choose 2}} x^k =(x;q)_n, \label{eq:qbinomial}
\end{align}
we have
\begin{align*}
\sum_{k=0}^{n}(-1)^{k}{n\brack k}q^{{k\choose 2}} q^{-ik}
=0 \quad\text{for}\ 0\leqslant i\leqslant n-1.  
\end{align*}
It follows that the expression \eqref{eq:last-lem-4} is equal to $0$. Namely, the identity \eqref{eq:last-lem-3} holds.
\end{itemize}
Hence, we see that all the $m+n$ roots of $L_{m,n}(x)$ are the same as those of the right-hand side of \eqref{eq:last-lem}.  \qed

\begin{lem}
Let $n$ and $h$ be positive integers and let $m$ be a non-negative integer with $h\leqslant n-m$. Then
\begin{align}
&\hskip -2mm\sum_{j=0}^{m}\sum_{k=m+h}^n \frac{(q^{-n};q)_{j}  (q^{-n};q)_{k}(a;q)_{j} (a;q)_{k} (1-q^{k-j}) q^{j+k+jh}  }
{(q;q)_{j} (q;q)_{k}(c;q)_{j} (c;q)_{k}}
{k-m-1\brack h-1}{m+h-j-1\brack h-1} \notag\\
&=\frac{(q;q)_n (a;q)_{m+h}(c/a;q)_{n-h}a^{n-h}q^{\frac{m^2+m-h^2+h}{2}-mn}}
{(-1)^{m-h}(q;q)_{m} (c;q)_{m} (c;q)_{n}(q;q)_{h-1}(q;q)_{n-m-h}}. \label{eq:am-2}
\end{align}
\end{lem}
\pf It is easy to see that ${k-m-1\brack h-1}=0$ for $m+1\leqslant k<m+h$.
Therefore, the left-hand side of \eqref{eq:am-2} remains unchanged when we replace $\sum_{k=m+h}^n$ by $\sum_{k=m+1}^n$.
Moreover,
\begin{align*}
{k-m-1\brack h-1}{m+h-j-1\brack h-1}
=\frac{(q^{j-m-h+1};q)_{h-1}(q^{k-m-h+1};q)_{h-1}q^{(m-j)(h-1)+{h\choose 2}} }{(-1)^{h-1}(q;q)_{h-1}^2}.
\end{align*}
The proof then follows from \eqref{eq:last-lem} and \eqref{eq:last-lem-2} with $x=a$.
\qed

The following result has been  proved  in \cite[(3.5)]{GZ2}.
\begin{lem}[\cite{GZ2}] Let $n$ be a positive integer. Then
\begin{align}
(x;q)_n+(a/x;q)_n &=(x;q)_n(a/x;q)_n+(a;q)_n
+\sum_{k=1}^{n-1} (x;q)_k (a/x;q)_k  B_{n,k}(a),
 \label{eq:lem-important2}
\end{align}
where
$$
B_{n,k}(a):=
(1-q^n)\sum_{h=1}^{n-k} (-1)^h {n-k-1\brack h-1}{k+h-1\brack h-1}\frac{ q^{{h\choose 2}+kh}a^h}{1-q^h}.
$$
\end{lem}
\section{Proof of Theorem {\rm\ref{thm:1.4}}  \label{sec:thm:1.4}}
As $$\left(\sum_{k=0}^n a_k\right)\left(\sum_{j=0}^n b_j\right)=\sum_{k=0}^n a_kb_k+\sum_{0\leqslant j<k\leqslant n}(a_kb_j+a_jb_k),$$
the left-hand side of \eqref{eq:thm1.4} is equal to
\begin{align}
&\sum_{k=0}^{n} \frac{(q^{-n};q)_k^2 (a;q)_k^2 q^{2k} }{(q;q)_{k}^2 (c;q)_{k}^2}(x;q)_k (c/x;q)_k \notag\\
&\quad{}+\sum_{0\leqslant j<k\leqslant n}
\frac{(q^{-n};q)_j (q^{-n};q)_k (a;q)_j (a;q)_k q^{j+k} \big((x;q)_j (c/x;q)_k+(x;q)_k (c/x;q)_j \big)}
{(q;q)_{j} (q;q)_{k}(c;q)_{j} (c;q)_{k} }.  \label{eq:thm-pf-1}
\end{align}
For $0\leqslant j<k$, from \eqref{eq:lem-important2} we deduce that
\begin{align*}
&\hskip -2mm (x;q)_j (c/x;q)_k+(x;q)_k (c/x;q)_j \\
&=(x;q)_j(c/x;q)_j\big((xq^j;q)_{k-j}+(cq^{j}/x;q)_{k-j}\big)  \\
&=(x;q)_k(c/x;q)_k+(x;q)_j(c/x;q)_j(cq^{2j};q)_{k-j} +\sum_{i=1}^{k-j-1} (x;q)_{j+i} (c/x;q)_{j+i} B_{k-j,i}(cq^{2j})\\
&=(x;q)_k(c/x;q)_k+(x;q)_j(c/x;q)_j +\sum_{i=0}^{k-j-1} (x;q)_{j+i} (c/x;q)_{j+i} B_{k-j,i}(cq^{2j}),
\end{align*}
where we have used the $q$-binomial theorem~\eqref{eq:qbinomial} in the last step:
$$
(cq^{2j};q)_{k-j}=1+\sum_{h=1}^{k-j}(-1)^h{k-j\brack h}q^{{h\choose 2}+2jh}c^h.
$$
It follows that \eqref{eq:thm-pf-1} can be written as
$
\sum_{m=0}^n \alpha_m (x;q)_m (c/x;q)_m,
$
where
\begin{align}
\alpha_m&=\sum_{j=0}^n\frac{(q^{-n};q)_j (q^{-n};q)_m (a;q)_j (a;q)_m  q^{j+m}}
{(q;q)_{j} (q;q)_{m}(c;q)_{j} (c;q)_{m} }  \notag\\
&\quad{}+\sum_{j=0}^{m}\sum_{k=m+1}^n \frac{(q^{-n};q)_j (q^{-n};q)_k (a;q)_j (a;q)_k  q^{j+k}  }
{(q;q)_{j} (q;q)_{k} (c;q)_{j}(c;q)_{k}}B_{k-j, m-j}(cq^{2j}).
\label{eq:am-0}
\end{align}
By  the $q$-Chu-Vandermonde summation formula~\cite[Appendix (II.6)]{GR}:
\begin{align}\label{q-chu}
_{2}\phi_{1}\left[\begin{array}{c}
a,q^{-n}\\
c
\end{array};q, q
\right]
=\frac{(c/a;q)_n}{(c;q)_n}a^n,
\end{align}
we have
\begin{align}
&\sum_{j=0}^n\frac{(q^{-n};q)_j (q^{-n};q)_m (a;q)_j (a;q)_m  q^{j+m}}
{(q;q)_{j} (q;q)_{m}(c;q)_{j} (c;q)_{m} }\nonumber\\
&=(-1)^m \frac{(q;q)_{n} (a;q)_m (c/a;q)_n q^{\frac{m^2+m}{2}-mn} } {(q;q)_m (c;q)_m (c;q)_n (q;q)_{n-m}}a^n.
\label{eq:am-1}
\end{align}
Substituting \eqref{eq:am-1} and \eqref{eq:am-2} into \eqref{eq:am-0}, we obtain
\begin{align}
\alpha_m&=\frac{(-1)^{m}(q;q)_n q^{\frac{m^2+m}{2}-mn}}{(q;q)_{m}(c;q)_{m}(c;q)_{n}}
\sum_{h=0}^{n-m}\frac{(a;q)_{m+h}(c/a;q)_{n-h}a^{n-h} q^{mh} c^h }
{(q;q)_{h}(q;q)_{n-m-h}}. \label{eq:fin-1}
\end{align}
The last sum can be summed again by the $q$-Chu-Vandermonde formula~\eqref{q-chu} and
is equal to
\begin{align}
\frac{(a;q)_{m} (c/a;q)_m (cq^{2m};q)_{n-m} a^n}{(q;q)_{n-m}}. \label{eq:fin-2}
\end{align}
It follows from \eqref{eq:fin-1} and \eqref{eq:fin-2} that $\alpha_m$ is just the coefficient of $(x;q)_m (c/x;q)_m$ on the right-hand side of
\eqref{eq:thm1.4}.

\medskip
\noindent{\it Remark.}
Letting $c=q^{2s+1}$, $a=aq^{s}$, $x=xq^{s}$, and replacing $n$ by $n-s$ in \eqref{eq:thm1.4} ($0\leqslant s\leqslant n$), we get the following result:
\begin{align}
&\hskip -3mm \left(\sum_{k=s}^{n} \frac{(q^{-n};q)_k (a;q)_k(x;q)_k q^k }{(q;q)_{k-s} (q;q)_{k+s}}\right)
  \left(\sum_{k=s}^{n} \frac{(q^{-n};q)_k (a;q)_k(q/x;q)_k q^k }{(q;q)_{k-s} (q;q)_{k+s}}\right) \nonumber \\
&=\frac{(q^{-n};q)_s(a;q)_s a^{n-s}q^{(n+1)s-s^2}}{(q^{n+1};q)_s (q/a;q)_s}
\sum_{k=s}^{n} \frac{(q^{-n};q)_k (q^{n+1};q)_k (a;q)_k (q/a;q)_k (x;q)_k (q/x;q)_k  q^{k}}{(q;q)_{k-s}(q;q)_{k+s} (q;q)_{2k}}. \label{eq:general-s}
\end{align}
It is clear that the $a=-q^{-n}$ case of \eqref{eq:general-s} reduces to \eqref{eq:important-thm-0}.

\section{Proof of Theorem~\ref{thm:last}\label{Pr:thm:last}}
We need a special case of Theorem~{\rm\ref{thm:1.4}.
Letting $a=-x$ in \eqref{eq:thm1.4}, we are led to
\begin{align}
&\hskip -3mm \left(\sum_{k=0}^{n} \frac{(q^{-n};q)_k (x^2;q^2)_k q^k }{(q;q)_{k} (c;q)_{k}}\right)
  \left(\sum_{k=0}^{n} \frac{(q^{-n};q)_k (-x;q)_k(c/x;q)_k q^k }{(q;q)_{k} (c;q)_{k}}\right) \nonumber \\
&=(-x)^n
\sum_{k=0}^{n} \frac{(q^{-n};q)_k (cq^{n};q)_k (x^2;q^2)_k (c^2/x^2;q^2)_k  q^{k}}{(q;q)_{k}(c;q)_{k} (c;q)_{2k}}. \label{eq:special-1}
\end{align}
We also need the following result.
\begin{lem}Let $n$ be a non-negative integer. Then
\begin{align}
\sum_{k=0}^{n} \frac{(q^{-n};q)_k (x;q)_k(y;q)_k q^k }{(q;q)_{k} (c;q)_{k}}
=\sum_{k=0}^{n} (-1)^k\frac{(q^{-n};q)_k (x;q)_k (c/y;q)_k x^{n-k }y^k q^{nk-{k\choose 2}} }{(q;q)_{k} (c;q)_{k}}. \label{eq:special-222}
\end{align}
\end{lem}
\pf
This follows from combining Jackson's two  transformations of terminating
$_2\phi_1$ series~\cite[Appendix (III.7) and (III.8)]{GR}:
\begin{align*}
{}_{3}\phi_{2}\left[\begin{array}{c}
q^{-n}, b, bzq^{-n}/c\\
bq^{1-n}/c, 0
\end{array};q,q
\right]
=b^n{}_{3}\phi_{1}\left[\begin{array}{c}
q^{-n}, b, q/z\\
bq^{1-n}/c
\end{array};q,\frac{z}{c}
\right]
\end{align*}
with
$b\to x,\;c\to xq^{1-n}/c$ and $z\to qy/c$. \qed

Replacing $x$ and $y$ by $-x$ and $c/x$ respectively in \eqref{eq:special-222}, we are led to
\begin{align}
\sum_{k=0}^{n} \frac{(q^{-n};q)_k (-x;q)_k(c/x;q)_k q^k }{(q;q)_{k} (c;q)_{k}}
=(-x)^n
\sum_{k=0}^{n} \frac{(q^{-n};q)_k (x^2;q^2)_k c^k q^{nk-{k\choose 2}} }{(q;q)_{k} (c;q)_{k} x^{2k}}. \label{eq:special-2}
\end{align}
Combining \eqref{eq:special-1} and \eqref{eq:special-2} (also $x\to\sqrt{x}$), we obtain
Theorem~\ref{thm:last}.

\medskip
\noindent{\it Remark.} Letting $c=q^{2s+1}$, $x\to xq^{2s}$, and replacing $n$ by $n-s$ ($0\leqslant s\leqslant n$)
in  \eqref{eq:special-3}, we get the following identity
\begin{align*}
&\hskip -3mm\left( \sum_{k=s}^{n} \frac{(q^{-n};q)_k (x;q^2)_k q^k }{(q;q)_{k-s} (q;q)_{k+s}} \right)
  \left( \sum_{k=s}^{n} \frac{(q^{-n};q)_k (x;q^2)_k q^{(n+1)k-{k\choose 2}}} {(q;q)_{k-s} (q;q)_{k+s} x^k} \right) \nonumber \\
&=\frac{(-1)^s(q;q)_n^2 (x;q^2)_s q^s}{(q;q)_{n-s}(q;q)_{n+s}(q^2/x;q^2)_s x^s}
\sum_{k=s}^{n} \frac{(q^{-n};q)_k (q^{n+1};q)_k (x;q^2)_k (q^2/x;q^2)_k  q^{k}}{(q;q)_{k-s}(q;q)_{k+s} (q;q)_{2k}},
\end{align*}
which was originally conjectured in a preliminary  version (arXiv:1408.0512v1) of \cite{GZ2}.

\section{Proof of Theorem~\ref{thm:1}\label{sec:proofoftheorem1} }
If we replace $n$ by $n-i$
in the $q$-Chu-Vandermonde summation
formula~\eqref{q-chu} with $a=q^{-m+i}$ and  $c=q^{i+1}$ then
$$
\sum_{k=i}^{m}q^{(m-k)(n-k)}{m\brack k}{n-i\brack k-i} ={n+m-i\brack n}.
$$
Hence, by the $q$-binomial theorem \eqref{eq:qbinomial}  we have
\begin{align*}
\sum_{k=0}^{m}q^{(m-k)(n-k)}{m\brack k}{n\brack k}(x;q)_k
&=\sum_{k=0}^{m}q^{(m-k)(n-k)}{m\brack k}{n\brack k}\sum_{i=0}^k (-1)^i {k\brack i}q^{i\choose 2}x^i  \\[5pt]
&=\sum_{i=0}^{m}(-1)^i {n\brack i} q^{i\choose 2}x^i \sum_{k=i}^{m}q^{(m-k)(n-k)}{m\brack k}{n-i\brack k-i}  \\[5pt]
&=\sum_{i=0}^{m} q^{i\choose 2}{n\brack i}{n+m-i\brack n} (-x)^i.
\end{align*}
When $x=-1$ and  $x=-q$, we obtain \eqref{D1} and \eqref{D2},  respectively.
Now, letting $c=q$ and $a=q^{-m}$ in \eqref{eq:thm1.4}, we get
\begin{align*}
&\hskip -2mm \left(\sum_{k=0}^{m}q^{(m-k)(n-k)}{m\brack k}{n\brack k}(x;q)_k\right)
\left(\sum_{k=0}^{m}q^{(m-k)(n-k)}{m\brack k}{n\brack k}(q/x;q)_k\right)  \nonumber\\[5pt]
&=\sum_{k=0}^{m}q^{(m-k)(n-k)}{n+k\brack 2k}{m\brack k}{m+k\brack k}(x;q)_k(q/x;q)_k.  
\end{align*}
This entails the identity \eqref{eq:thm:1}  by taking $x=-1$ or $-q$.

\section{Proof of Theorem \ref{thm:2}\label{sec:pf2}}
The following identity  can be easily proved by induction.
\begin{align}
\sum_{k=j}^{n-1}(1-q^{2k+1}){k+j\brack 2j}q^{-(j+1)k}=\frac{(1-q^{n})(1-q^{n-j})}{1-q^{j+1}}{n+j\brack 2j}q^{-(j+1)(n-1)}.\label{eq:qidentity}
\end{align}
By \eqref{eq:thm:1} and \eqref{eq:qidentity}, the left-hand side of \eqref{eq:cases} is equal to
\begin{align}
&\hskip -2mm \sum_{k=0}^{p-1}\frac{1-q^{2k+1}}{1-q}
\sum_{j=0}^{k}q^{j^2-mj-(j+1)k}{k+j\brack 2j}{m\brack j}{m+j\brack j}(-1;q)_j(-q;q)_j  \notag\\[5pt]
&=\sum_{j=0}^{p-1}q^{j^2-mj-(j+1)(p-1)}\frac{(1-q^{p})(1-q^{p-j})}{(1-q)(1-q^{j+1})}{p+j\brack 2j}{m\brack j}{m+j\brack j}(-1;q)_j(-q;q)_j. \label{eq:sum-thm2}
\end{align}
By \cite[Theorem 2.1]{Guo}, we know that
\begin{align*}
&\hskip -2mm\frac{(1-q^m)(1-q^{m+1})(1-q^{p-j})}{(1-q)(1-q^p)(1-q^{j+1})}{p+j\brack 2j}{m\brack j}{m+j\brack j} \\[5pt]
&=\frac{(1-q^{p-j})(1-q^{j+1})}{(1-q)(1-q^p)}{p+j\brack 2j}{m+1\brack j+1}{m+j\brack j+1}
\end{align*}
is a polynomial in $q$ with non-negative integer coefficients. Since $[p]=(1-q^p)/(1-q)$ is an irreducible polynomial in $q$ for any prime $p$
and $\gcd(q^m-1,q^n-1)=(q^{\gcd(m,n)}-1)$, we conclude that
$$
\frac{1-q^{p-j}}{1-q^{j+1}}{p+j\brack 2j}{m\brack j}{m+j\brack j}\equiv 0\pmod{[p]}\quad\text{if $m\not\equiv 0, -1\pmod{p},$}
$$
and so the right-hand side of \eqref{eq:sum-thm2} is congruent to $0$ modulo $[p]^2$ in this case.

On the other hand, if $m\equiv 0,-1\pmod{p}$, then
$$
\frac{1-q^{p-j}}{1-q^{j+1}}{m\brack j}{m+j\brack j}
\equiv 0\pmod{[p]}\quad\text{if $j=0,1,\ldots,p-2.$}
$$
Therefore, if $m\equiv 0\pmod{p}$, then the right-hand side of \eqref{eq:sum-thm2} is congruent to
\begin{align*}
q{2p-1\brack 2p-2}{m\brack p-1}{m+p-1\brack p-1}(-1;q)_{p-1}(-q;q)_{p-1}
&\equiv q\frac{1-q^{2p-1}}{1-q}.\frac{1-q^m}{1-q^{p-1}}(-q^{-1})\cdot \frac{2q}{1+q} \\[5pt]
&\equiv -\frac{2q(1-q^m)}{1-q^2} \\[5pt]
&\equiv \frac{1-q^{-2m}}{1-q^2}q \pmod{[p]^2},
\end{align*}
where we have used the congruence $(-q;q)_{p-1}\equiv 1\pmod{[p]}$ (see \cite[(1.6)]{Guo2015} or \cite{Pan});
while if $m\equiv -1\pmod{p}$, then the right-hand side of \eqref{eq:sum-thm2} is congruent to
\begin{align*}
\frac{1-q^{2p-1}}{1-q}.\frac{1-q^{m+1}}{1-q^{p-1}}\cdot \frac{2q}{1+q}
&\equiv \frac{2q(1-q^{m+1})}{1-q^2} \\[5pt]
&\equiv \frac{1-q^{2m+2}}{1-q^2}q \pmod{[p]^2}.
\end{align*}

\section{Proof of Theorem \ref{thm:3} \label{sec:pf3}}
Similarly to \eqref{eq:sum-thm2}, the left-hand side of \eqref{eq:thm3-1} is equal to
\begin{align*}
\sum_{j=0}^{n-1}q^{j^2-mj-(j+1)(n-1)}\frac{(1-q^{m})(1-q^{m+1})(1-q^{n-j})}{(1-q^2)(1-q^n)(1-q^{j+1})}{n+j\brack 2j}{m\brack j}{m+j\brack j}
(-1;q)_j(-q;q)_j.
\end{align*}
It is easy to see that
\begin{align*}
&\hskip -2mm\frac{(1-q^{m})(1-q^{m+1})(1-q^{n-j})}{(1-q^2)(1-q^n)(1-q^{j+1})}{n+j\brack 2j}{m\brack j}{m+j\brack j}(-q;q)_j \\[5pt]
&=\begin{cases}
\displaystyle\frac{(1-q^{m})(1-q^{m+1})}{(1-q^2)(1-q)}&\text{if $j=0$,}\\[10pt]
\displaystyle\frac{(1-q^{n-j})(1-q^{j+1})}{(1-q)(1-q^n)}{n+j\brack 2j}{m+1\brack j+1}{m+j\brack j+1}(-q^2;q)_{j-1}
&\text{if $j\geqslant 1$}
\end{cases}
\end{align*}
is a polynomial in $q$ with non-negative integer coefficients by \cite[Theorem 2.1]{Guo}. We conclude that
\eqref{eq:thm3-1} is the desired Laurent polynomial in $q$.

Let
$$
S_n(x_0,\ldots,x_n;q)=\sum_{k=0}^{n}{n+k\brack 2k}{2k\brack k}q^{-nk} x_k.
$$
To prove that \eqref{eq:thm3-2} and \eqref{eq:thm3-3} also have the same properties, we first establish the following result.
\begin{lem}\label{lem:4.1}
Let $n$ and $r$ be positive integers. Then both
\begin{align*}
\sum_{k=0}^{n-1}\frac{1-q^{2k+1}}{1-q^n} S_k(x_0,\ldots,x_k)^r q^{-k}\ \text{and}\
\sum_{k=0}^{n-1}(-1)^{n-k-1}\frac{1-q^{2k+1}}{1-q^n} S_k(x_0,\ldots,x_k)^r q^{k\choose 2}
\end{align*}
are polynomials in $x_0,\ldots,x_{n-1},q$ and $q^{-1}$ with non-negative integer coefficients.
\end{lem}
\pf Recall the identity
\begin{align*}
{k+i\brack 2i}{2i\brack i}{k+j\brack 2j}{2j\brack j}=\sum_{s=i}^{i+j}{i+j\brack i}{j\brack s-i}{s\brack j}{k+s\brack 2s}{2s\brack s}q^{(i+j-s)(k-s)}.
\end{align*}
which can be proved using  the $q$-Pfaff-Saalsch\"utz identity (see \cite[Lemma 2.1]{Schmidt}).
It follows that
\begin{align}
S_k(x_0,\ldots,x_k)^r
&=\sum_{0\leqslant i_1,\ldots,i_r\leqslant k}\prod_{j=1}^{r}{k+i_j\brack 2i_j}{2i_j\brack i_j} q^{-ki_j}x_{i_j} \nonumber\\[5pt]
&=\sum_{0\leqslant i_1,\ldots,i_r\leqslant k}x_{i_1}\cdots x_{i_r}\sum_{s=i_1}^{i_1+\cdots+i_r}
P(i_1,\ldots,i_r,s){k+s\brack 2s}{2s\brack s} q^{-ks},  \label{eq:sk-power}
\end{align}
where $P(i_1,\ldots,i_r,s)$ is a Laurent polynomial in $q$ independent of $k$ with non-negative integer coefficients. Therefore, by
\eqref{eq:qidentity}, we see that
\begin{align*}
&\hskip -2mm\sum_{k=0}^{n-1}\frac{1-q^{2k+1}}{1-q^n} S_k(x_0,\ldots,x_k)^r q^{-k} \\[5pt]
&=\sum_{0\leqslant i_1,\ldots,i_r\leqslant n-1}x_{i_1}\cdots x_{i_r}\sum_{s=i_1}^{i_1+\cdots+i_r}
P(i_1,\ldots,i_r,s)\frac{1-q^{n-s}}{1-q^{s+1}}{n+s\brack 2s}{2s\brack s} q^{-(s+1)(n-1)}
\end{align*}
is a polynomial in $x_0,\ldots,x_{n-1},q$ and $q^{-1}$ with non-negative integer coefficients
since
$$\frac{1-q^{n-s}}{1-q^{s+1}}{n+s\brack 2s}{2s\brack s}={n+s\brack s}{n\brack s+1}.$$

Similarly, since
\begin{align*}
\sum_{k=s}^{n-1}(-1)^{n-k-1}\frac{1-q^{2k+1}}{1-q^n}{k+s\brack 2s}{2s\brack s}q^{{k\choose 2}-sk}
={n-1\brack s}{n+s\brack s}q^{{n\choose 2}-sn},
\end{align*}
we deduce  from \eqref{eq:sk-power} that\begin{align*}
&\hskip -2mm\sum_{k=0}^{n-1}(-1)^{n-k-1}\frac{1-q^{2k+1}}{1-q^n} S_k(x_0,\ldots,x_k)^r q^{k\choose 2} \\[5pt]
&=\sum_{0\leqslant i_1,\ldots,i_r\leqslant n-1}x_{i_1}\cdots x_{i_r}\sum_{s=i_1}^{i_1+\cdots+i_r}
P(i_1,\ldots,i_r,s){n-1\brack s}{n+s\brack s}q^{{n\choose 2}-sn}
\end{align*}
is a polynomial in $x_0,\ldots,x_{n-1},q$ and $q^{-1}$ with non-negative integer coefficients.
\qed

For $k=0,\ldots,n-1$, let
$$
x_k={m+k\brack 2k}(-1;q)_k(-q;q)_k q^{k^2-mk}.
$$
Then the identity \eqref{eq:thm:1} may be rewritten as
\begin{align*}
D_q(m,n)D_{q^{-1}}(m,n)=\sum_{k=0}^{n}{n+k\brack 2k}{2k\brack k}q^{-nk}x_k.
\end{align*}
It is clear that $x_0,\ldots,x_{n-1}$ are Laurent polynomials in $q$ with non-negative integer coefficients.
By Lemma \ref{lem:4.1}, so are the expressions \eqref{eq:thm3-2} and \eqref{eq:thm3-3}.

\vskip 5mm
\noindent{\bf Acknowledgments.} The authors would like to thank the referees and the editor for helpful comments on a previous version of this paper.
The second author was partially
sponsored by the National Natural Science Foundation of China (grant 11371144),
the Natural Science Foundation of Jiangsu Province (grant BK20161304),
and the Qing Lan Project of Education Committee of Jiangsu Province.

\end{document}